\documentstyle[amssymb,amsfonts]{amsart}

\theoremstyle{plain}
  \newtheorem{theorem}[subsection]{Theorem}
  \newtheorem*{theorem*}{Theorem}

  \newtheorem{lemma}[subsection]{Lemma}
  \newtheorem{corollary}[subsection]{Corollary}

\theoremstyle{remark}

\theoremstyle{definition}

\begin{document}

\title[On a question of Erd\H os and Ulam]
{On a question of Erd\H os and Ulam}
\author{Jozsef Solymosi \and Frank de Zeeuw}
\address{Department of Mathematics, UBC, 1984 Mathematics Road, Vancouver, BC, Canada V6T 1Z2}
\email{\{fdezeeuw,solymosi\}@@math.ubc.ca}
\thanks{The first author was supported by NSERC and OTKA grants and by Sloan Research Fellowship}

\begin{abstract}
Ulam asked in 1945 if there is an everywhere dense \emph{rational
set}, i.e. a point set in the plane with all its pairwise
distances rational. Erd\H os conjectured that if a set $S$ has a
dense rational subset, then $S$ should be very special. The only
known types of examples of sets with dense (or even just infinite)
rational subsets are lines and circles. In this paper we prove
Erd\H os' conjecture for algebraic curves, by showing that no
irreducible algebraic curve other than a line or a circle contains
an infinite rational set.
\end{abstract}

\maketitle

\section{Introduction}
We define a {\it rational set} to be a set $S\subset {\mathbb
R}^2$ such that the distance between any two elements is a
rational number. We are interested in the existence of infinite
rational distance sets on algebraic curves.

On any line, one can easily find an infinite rational set that is
in fact dense. It is an easy exercise to find an everywhere dense
rational subset of the unit circle. On the other hand, it is not
known if there is a rational set with 8 points {\it in general
position}, i.e. no 3 on a line, no 4 on a circle. In 1945, Anning
and Erd\H os \cite{AE} proved that any infinite \emph{integral}
set, i.e. where all distances are integers, must be contained in a
line. Problems related to rational and integral sets became one of
Erd\H os' favorite subjects in combinatorial geometry \cite{E3}
\cite{E4} \cite{E5} \cite{E6} \cite{E7} \cite{E8}.

In 1945, when Ulam heard Erd\H os' simple proof \cite{E1} of his
theorem with Anning, he said that he believed there is no
everywhere dense rational set in the plane, see Problem III.5. in
\cite{UL} and also \cite{E2}. Erd\H os conjectured that an
infinite rational set must be very restricted, but that it was
probably a very deep problem \cite{E2}\cite{E3}. Not much progress
has been made on Ulam's question. There were attempts to find
rational sets on parabolas \cite{Ca} \cite{Ch}, and there were
some results on integral sets, in particular bounds were found on
the diameter of integral sets \cite{HKM} \cite{So}. Very recently
Kreisel and Kurz \cite{KK} found an integral set with 7 points in
general position.

In this paper, we prove that lines and circles are the only
irreducible algebraic curves that contain infinite rational sets.
Our main tool is Faltings' Theorem \cite{Fal}. We will also show
that if a rational set $S$ has infinitely many points on a line or
on a circle, then all but 4 resp. 3 points of $S$ are on the line
or on the circle. This answers questions of Guy, Problem D20 in
\cite{Guy}, and Pach, Section 5.11 in \cite{BMP}.

\section{Main result}
Our main result is the following.
\begin{theorem}\label{main}
Every rational set of the plane has only finitely many points in
common with an algebraic curve defined over $\mathbb{R}$, unless
the curve has a component which is a line or a circle.
\end{theorem}
The two special cases, line and circle, are treated in more detail
in the next theorem.
\begin{theorem}\label{minor}
If a rational set $S$ has infinitely many points on a line or on a
circle, then all but 4 resp. 3 points of $S$ are on the line or on
the circle.
\end{theorem}
Note that there are infinite rational sets with all but 4 points
on a line, and there are infinite rational sets with all but 3
points on a circle. The circle case follows from the line case by
applying an inversion with rational radius and center one of the 4
points not on the line. A construction of Huff \cite{Hu}\cite{Pe}
gives an infinite rational set with all but 4 points on a line.

We can formulate our Theorem \ref{main} in a different way by
using the term \emph{curve-general position}: we say that a point
set $S$ of ${\mathbb R}^2$ is in curve-general position if no
algebraic curve of degree $d$ contains more than $d(d+3)/2$ points
of $S$. Note that $d(d+3)/2$ is the number of points in general
position that determine a unique curve of degree $d$.

\begin{corollary}
If $S$ is an infinite rational set in general position, then there
is an infinite $S'\subset S$ such that $S'$ is in curve-general
position.
\end{corollary}

{\it Proof:} Let $S_5$ consist of any five points in $S$, and let
$T_5$ be the set of finitely many points on the unique conic
through those 5 points. Continue recursively; at step $n$, add a
point from $S\backslash T_{n-1}$ to $S_{n-1}$ to get $S_n$. For
each $d$ such that $d(d+3)/2\leq n$, let $T_n$ be the union of
$T_{n-1}$ and the set of points of $S$ that are on a curve of
degree $d$ through any $d(d+3)/2$ points in $S_{n}$. Since each
$T_n$ is finite, we can always add another point. Then the
infinite union of the sets $S_n$ is an infinite subset of $S$ with
the required property.

\section{Proof of Theorem \ref{main}}

\subsection{General Approach}
We will use the following theorem of Faltings \cite{Fal}.
\begin{theorem*}[Faltings]
A curve of genus $\geq 2$, defined over a number field, contains
only finitely many rational points.
\end{theorem*}
\noindent In this paper by \emph{curve} (defined over a field
$K\subset \mathbb{R}$) we usually mean the zero set in ${\mathbb
R}^2$ of a polynomial in two variables with coefficients from $K$.
But when we consider the genus of a curve, we are actually talking
about the projective variety defined by the polynomial. For
definitions, see \cite{Sil}.

First suppose we have an infinite rational set $S$ contained in a
curve $C$ of genus $\geq 2$, defined over $\mathbb{R}$. We can
move two points in $S$ to $(0,0)$ and $(0,1)$, so that by Lemma
\ref{lemma2} below the elements of $S$ are of the form $(r_1,
r_2\sqrt{k})$. Then by the remark after Lemma \ref{lemma2}, the
curve is defined over ${\mathbb Q}(\sqrt{k})$. By Faltings'
theorem, $S$ must be finite.

Below we will show that if we have an infinite rational set $S$ on
a curve $C_1$ of genus 0 or 1, then all but finitely many of the
points in $S$ will in fact give points on a curve $C_2$ in
${\mathbb R}^3$ of genus $\geq 2$. More precisely, assuming
$(0,0)$ and $(0,1)$ are in $S$, a point $(r_1,r_2\sqrt{k})$ will
give a point $(r_1,r_2\sqrt{k},r_3)$ on a curve $C_2$, with all
the $r_i$ rational. Again we conclude by Faltings' theorem that
the original set $S$ could not have been infinite.

\subsection{Two lemmata}
Rationality of distances in ${\mathbb R}^2$ is clearly preserved
by translations, rotations and uniform scaling, $(x,y)\mapsto
(\lambda x,\lambda y)$ with $\lambda\in \mathbb{Q}$. More
surprisingly, rational sets are preserved under certain central
inversions. This will be an important tool in our proof below.

\begin{lemma}\label{lemma1}
If we apply inversion to a rational set $S$, with center a point
$x\in S$ and rational radius, then the image of $S\backslash\{x\}$
is a rational set.
\end{lemma}

{\it Proof:} We may assume the center to be the origin and the
radius to be 1. The properties of inversion are most easily seen
in complex notation, where the map is $z\mapsto 1/z$. Suppose we
have two points $z_1$, $z_2$ with rational distances $|z_1|$,
$|z_2|$ from the origin, and with $|z_2-z_1|$ rational. Then
$$\left|\frac{1}{z_1} - \frac{1}{z_2}\right| = \left|
\frac{z_2-z_1}{z_1 z_2}\right| = \frac{|z_2-z_1|}{|z_1||z_2|}$$ is
also rational.\medskip

A priori, points in a rational set could take any form. However,
after moving two of the points to two fixed rational points by
translating, rotating, and scaling, the points are in fact almost
rational points. The following simple lemma is well-known among
researchers working with integer sets. As far as we know, it was
proved first by Kemnitz \cite{Ke}.

\begin{lemma}\label{lemma2}
For any rational set $S$ there is a square-free integer $k$ such
that if a similarity transformation $T$ transforms two points of
$S$ into $(0,0)$ and $(1,0)$ then any point in $T(S)$  is of the
form
$$(r_1, r_2\sqrt{k} ),~~~r_1,r_2\in\mathbb{Q}.$$
\end{lemma}

Note that this implies that any curve of degree $d$ containing at
least $d(d+3)/2$ points from $T(S)$ is defined over ${\mathbb
Q}(\sqrt{k})$.

\subsection{Curves of genus 1}
Let $C_1:f(x,y)=0$ be an irreducible algebraic curve of genus
$g_1=1$ and degree $d\geq3$. Suppose that there is an infinite set
$S$ on $C_1$ with pairwise rational distances. Assume that the
points $O=(0,0)$ and $(1,0)$ are on $C_1$ and in $S$, and that $O$
is not a singularity of $C_1$. Below we will be allowed to make
any other assumptions on $C_1$ that we can achieve by translating,
rotating or scaling it, as long as we also satisfy the assumptions
above. In particular, we can use any of the infinitely many rotations
about the origin that put another point of $S$ on the $x$-axis.

We wish to show that the intersection curve $C_2$ of the surfaces
\begin{align*}
X: & ~f(x,y)=0,\\
Y: & ~x^2+y^2=z^2,
\end{align*}
has genus $g_2 \geq 2$.

Consider $C_1$ as a curve in the $z=0$ plane, and define the map
$\pi:C_2\to C_1$ by $(x,y,z)\mapsto (x,y)$, i.e. the restriction
to $C_2$ of the vertical projection from the cone $Y$ to the $z=0$
plane. The preimage of a point $(x,y)$ consists of the two points
$(x,y,\pm\sqrt{x^2+y^2})$, except when $x^2+y^2=0$, which in
${\mathbb C}^2$ happens on the two lines $x+i y=0$ and $x-iy=0$.
Then we can determine (or at least bound from below) the genus of
$C_2$ using the Riemann-Hurwitz formula \cite{Sil} applied to
$\pi$,
$$2g_2-2 \geq \deg\pi\cdot(2g_1-2) + \sum_{P\in C_2} (e_P-1).$$
This is usually stated with equality for smooth curves, but we are
allowing $C_1$ and $C_2$ to have singularities. To justify our use
of it, observe that the map $\pi$ corresponds to a map
$\widetilde{\pi}:\widetilde{C}_1\to \widetilde{C}_2$ between the
normalizations of the curves, for which Riemann-Hurwitz holds. The
normalizations have the same genera as the original curves, and
$\widetilde{\pi}$ has the same degree. Furthermore a ramification
point of $\pi$ away from any singularities gives a ramification
point of $\widetilde{\pi}$. It is enough for our purposes to have
this inequality, but there could be more ramification points for
$\widetilde{\pi}$, above where the singularities used to be.

Applying this formula with $g_1=1$, $d=2$, we have
$$g_2 \geq 1+\frac{1}{2}\sum_{P\in C_2} (e_P-1),$$
so to get $g_2\geq 2$, we only need to show that $\pi$ has some
ramification point.

The potential ramification points are above the intersection
points of $C_1$ with the lines $x\pm iy=0$, of which there are
$2d$ by B\'ezout's theorem, counting with multiplicities. Such an
intersection point $P$ can only fail to have a ramification point
above it if the curve has a singularity at $P$, or if the curve is
tangent to the line there. We will show that there are only
finitely many lines through the origin on which one of those two things 
happens. Then certainly one of the infinitely many rotations of
$C_1$ that we allowed above will give an intersection point of
$C_1$ with $x\pm iy=0$ that has a ramification point above it.

The intersection of a line $y=a x$ with $f(x,y)=0$ is given by
$p_a(x)=f(x,a x)=0$, and if the point of intersection is a
singularity or a point of tangency, then $p_a(x)$ has a multiple
root. We can detect such multiple roots by taking the discriminant
of $p_a(x)$, which will be a polynomial in $a$ that vanishes if
and only if $p_a(x)$ has a multiple root. Hence for all but
finitely many values of $a$ the line $y=ax$ has $d$ simple
intersection points with $f(x,y)=0$. So indeed there is an allowed
rotation after which $\pi$ is certain to have a ramification
point.

\subsection{Curves of genus 0, $d\geq 4$}
Let $C_1:f(x,y)=0$ be an irreducible algebraic curve of genus
$g_1=0$, and again assume that it passes through the origin, but
does not have a singularity there. Then Riemann-Hurwitz with the
same map $\pi$ as above gives
$$g_2\geq -1 +\frac{1}{2}\sum_{P\in C_2} (e_P-1),$$
so to get $g_2\geq 2$ we need to show that there are at least 5
ramification points. As above, we can ensure that the lines $x\pm
iy$ each have $d$ simple points of intersection. Discounting the
intersection point of the two lines, this gives $2d-2$
ramification points. Hence if the degree of $f$ is $d\geq 4$ we
are done.

\subsection{Curves of genus 0, $d=2,3$.}

Let $d=3$ and assume $f(x,y)=0$ is not a line or a circle.
Consider applying inversion with the origin as center to the
curve. This is a birational transformation, so does not change the
genus. Therefore, when inversion increases the degree of $f$ to 
above 4, we are done.

Algebraically, inversion in the circle around the origin with radius 1 is given
by
$$(x,y)\mapsto \left(\frac{x}{x^2+y^2}, \frac{y}{x^2+y^2}\right),$$
and since this map is its own inverse, the curve $f(x,y)=0$ is sent to the curve
$$C_3:(x^2+y^2)^k\cdot f\left(\frac{x}{x^2+y^2}, \frac{y}{x^2+y^2}\right)=0,$$
where $k\leq d$ is the lowest integer that makes this a polynomial. This curve
is irreducible if and only if the original curve is irreducible.
Since $f$ does not have a singularity at the origin, it has a
linear term $a x+b y$ with $a$, $b$ not both zero. After inversion
this gives a highest degree term
$$(ax+by) (x^2+y^2)^{k - 1}.$$
In our situation, $d=3$, so if $k=3$, the curve $C_3$ has degree $2k-1=5$, and we are done.

The only other possibility is that $k=2$, which happens if $x^2+y^2$ divides the leading terms
of $f$. We will treat these cases in a completely different way. 

If $d=2$, then applying inversion will give a curve of degree 3, unless its leading terms are $x^2+y^2$, which
exactly means that it is a circle! So we treat this case by reducing it to the $d=3$ case. 

Since $f$ has degree 3 and genus 0, it must have a singularity. The singularity need not be in 
our rational set, but it is always a rational point, so we can move it to the origin, while maintaining
the almost-rational form of the points in our rational set. Then $f$ must have the form
$$(ax+by)(x^2+y^2) + cx^2+dy^2 + exy.$$
Note that this is exactly what we get if we apply inversion to a quadratic that is not a circle and goes through the origin.

In fact, we can ensure that $(1,0)$ is on the curve again, so that $a +c = 0$. Then if we divide by $c$, $f$ is of the form
$$(-x+by)(x^2+y^2) + x^2+dy^2 + exy.$$
We can parametrize this curve using lines $x = ty$, giving the parametrization
$$y(t) = \frac{t^2+et+d}{(t-b)(t^2+1)} =: \frac{p(t)}{q(t)},\hspace{30pt}x(t) = t\cdot y(t).$$
If we let $t_j$ be a value of $t$ that gives one of the points from our rational distance set, it follows that for infinitely many $t$,
$$\left(y(t) - y(t_j)\right)^2 + \left(x(t) - y(t_j)\right)^2 = \left( \frac{p(t)}{q(t)} - \frac{p(t_j)}{q(t_j)} \right)^2 + \left( t\cdot\frac{p(t)}{q(t)} - t_j\cdot\frac{p(t_j)}{q(t_j)} \right)^2$$
is a square. Then we can multiply by $q(t)^2q(t_j)^2$ to get infinitely many squares of the form
$$\left(p(t)q(t_j) - p(t_j)q(t)\right)^2 + \left(tp(t)q(t_j) - t_jp(t_j)q(t)\right)^2.$$
This polynomial has degree 6 in $t$. It has a factor $(t-t_j)^2$, and a factor $t^2+1$, since taking $t = \pm i$ gives (using $q(\pm i) = 0$)
$$\left(p(\pm i)q(t_j)\right)^2 + \left(\pm i \cdot p(\pm i)q(t_j)\right)^2 = 0.$$
Factoring these out, we get a quadratic polynomial $Q_j(t)$ in $t$. Its leading coefficient is
$$(t_j^2+1)((d^2+b^2)t_j^2+ 2(b^2e+ db - d^2b) t_j + b^2 e^2 + b^2 d^2 + d^2 + 2ebd),$$
and its constant term is
$$(t_j^2+1)((1+ (e+b)^2)t_j^2 + 2(bd - b +de) t_j+ d^2 + b^2).$$
These polynomials in $t_j$ are not identically zero (if $b$ and $d$ were both $0$, then $f$ would be reducible), hence we can pick $t_j$ so that they are not zero. Then in turn $Q_j(t)$ is a proper quadratic polynomial, and since it is essentially a distance function in the real plane, it cannot have real roots, so it has two distinct imaginary roots.

Therefore our infinite rational set gives infinitely many solutions to equations
$$z_j^2 = (t^2+1) \cdot Q_j(t).$$
Multiplying three of these together, and dividing out $(t^2+1)^2$, we get infinitely many solutions to
$$z^2 = (t^2+1)  Q_1(t)Q_2(t)Q_3(t).$$
If there are no multiple roots on the right, then this is a hyperelliptic curve of degree $8$, so it has genus $3$, hence cannot have infinitely many solutions, a contradiction.

The one thing we need to check is that we can choose the $t_j$ so that the $Q_j$ do not have roots in common. We need some notation: write
$$Q_j(t) = c_2(t_j) t^2 + c_1(t_j) t +c_0(t_j),$$
where
\begin{align*}
 c_2(t_j) & = (1+ (e+b)^2)t_j^2 + 2(bd +de- b ) t_j+ d^2 + b^2\\
 c_1(t_j) & = 2(bd+de - b)t_j^2 +2(b^2+d^2 - bed - bd - be - d)t_j + 2(bd + b^2e -bd^2)\\
 c_0(t_j) & = (d^2+b^2)t_j^2+ 2(b^2e+ db - d^2b) t_j + b^2 e^2 + b^2 d^2 + d^2 + 2ebd.
\end{align*}
Suppose that for infinitely many $t_j$ the polynomial $Q_j(t)$ has the same roots $x_1$ and $x_2$. Then for each of those $t_j$ we have 
$$c_1(t_j) = -(x_1+x_2)\cdot c_2(t_j),~~~c_0(t_j) = x_1\cdot x_2\cdot c_2(t_j).$$
If we look at the coefficients of the $t_j$ terms in these equations, we see that 
\begin{align*}
 -x_1-x_2 & = \frac{2(b^2+d^2 - bed - bd - be - d)}{2(bd - b +de)} = - b - \frac{be+d-d^2}{bd+de - b}\\
x_1\cdot x_2 & = \frac{2(b^2e+ db - d^2b)}{2(bd +de- b )} = b\cdot \frac{be+d-d^2}{bd+de - b}.
\end{align*}
Here we can read off that the roots are $x_1 =b$ and $x_2 =  \frac{be+d-d^2}{bd+de - b}$, which is a contradiction, since the roots had to be imaginary.

\section{Proof of Theorem \ref{minor}}
We will prove that if a rational set has infinitely many points on
a line, then it can have at most 4 points off the line.  The
corresponding statement for 3  points off a circle then follows by
applying an inversion. More precisely, suppose we have a rational
set $S$ with infinitely many points on a circle $C$ and at least 4
points off that circle. Assume that the origin is one of the
points in $S\cap C$, and apply inversion with the origin as
center, and with some rational radius. That turns $C$ into a line
$L$, and we get a rational set with infinitely many points on $L$,
and 4 other points. Moreover, the new origin can be added to $S$,
so that we get 5 points off the line, contradicting what we will
prove below. To see that the new origin has rational distance to
all points in $S$, observe that in complex notation the distances
$|z|$ to the old origin were rational for all $z\in S$, and that
the distances to the new origin are $1/|z|$.

To prove the statement for a line, our main tool will again be
Faltings' theorem, but now applied to a hyperelliptic curve
$$y^2 = \prod_{i=1}^6 (x-\alpha_i),$$
which has genus 2 if and only if the $\alpha_i$ are distinct.

Suppose we have a rational set $S$ with infinitely many points on
a line, say the $x$-axis, and 5 or more points off that line. Then
we can assume that 3 of those points are above the $x$-axis, and
that one of them is at $(0,1)$. Let the other two points be at
$(a_1,b_1)$ and $(a_2,b_2)$. Note that we are taking 3 points on
one side of the line, because we want to avoid having one point a
reflection of another. If we had, say, $(a_1,b_1)=(0,-1)$, the
argument below would break down.

Take a point $(x,0)$ of $S$ on the $x$-axis, with $x\neq 0,
a_1,a_2$. Then we have that
$$x^2+1,~~ (x-a_1)^2+b_1^2,~~\text{and}~~(x-a_2)^2+b_2^2$$
are rational squares, so that we get a rational point $(x,y)$ on
the curve
$$y^2 = (x^2+1)((x-a_1)^2+b_1^2)((x-a_2)^2+b_2^2).$$
To show that this is indeed a curve of genus 2, we need to show
that the right hand side does not have multiple roots. Since
$b_1,b_2\neq 0$, the 3 factors do not have real roots, hence each
has two distinct imaginary roots. Since each has the same
coefficient of the $x^2$ term, two of the quadratic polynomials
could only have common roots if they were identical. But the
coefficients of the $x$ term are respectively 0, $-2a_1$ and
$-2a_2$, so only the last two factors could be equal. But if
$a_1=a_2$, then $b_1\neq \pm b_2$, hence the constant terms would
be $a_1^2+b_1^2\neq a_2^2+b_2^2$. We conclude that the 3 factors
do not have common roots.

Therefore the curve has genus 2, and cannot contain infinitely
many rational points, contradicting the fact that $S$ has
infinitely many points on the line.

\section{Acknowledgements}
We thank Kalle Karu for the useful discussions. We are also indebted to an anonymous referee who noticed that the $d=3$ case in section 3.7 was not completely covered in the previous version of the paper.


\begin{thebibliography}{99}

\bibitem{AE} N.H.~Anning and P.~Erd\H os, \emph{Integral distances},
Bull. Amer. Math. Soc., (1945) 51, 598--600.

\bibitem{UL} S.M.~Ulam, \emph{A Collection of Mathematical Problems},
Interscience Tracts in Pure and Applied Mathematics, Number 8.
Interscience Publishers (1960) XIII, 150 p.

\bibitem{BMP} P.~Brass, W.~Moser, and J.~Pach, \emph{Research Problems
in Discrete Geometry.}
Springer, Berlin. 1st ed. (2005) XII, 499 p.

\bibitem{Ca} G.~Campbell, \emph{Points on $y = x^2$ at rational distance},
Math. Comp.,(2004) 73, 2093--2108.

\bibitem{Ch} A.~Choudhry, \emph{Points at Rational Distances on a
Parabola}, Rocky Mountain J. Math. Volume 36, Number 2 (2006),
413--424.

\bibitem{E1} P.~Erd\H os, \emph{Integral distances.},
Bull. Amer. Math. Soc., (1945) 51, 996.

\bibitem{E2} P.~Erd\H os, \emph{Ulam, the Man and the Mathematician}
J. Graph Theory 9 (1985) no. 4, 445--449. Also appears in {\it
Creation Math.} {\bf 19} (1986), 13--16.


\bibitem{E3} P.~Erd\H os,
\emph{Some combinatorial and metric problems in geometry},
Colloquia Mathematica Societatis J\'anos Bolyai 48 . Intuitive
Geometry, Si\'ofok, (1985) 167--177.

\bibitem{E4} P.~Erd\H os,
\emph{On Some Problems of Elementary and Combinatorial Geometry },
Annali di Matematica pura ed applicata, (IV), Vol. CIII, (1975),
99--108.

\bibitem{E5} P.~Erd\H os,
\emph{Verchu niakoy geometritchesky zadatchy},(On some geometric
problems, in Bulgarian), Fiz.-Mat. Spis. B\u ulgar. Akad. Nauk.
5(38) (1962), 205--212

\bibitem{E6} P.~Erd\H os,
\emph{Combinatorial problems in geometry}, Math. Chronicle 12
(1983), 35--54.

\bibitem{E7} P.~Erd\H os,
\emph{N\'eh\'any elemi geometriai probl\'em\'ar\'ol} (On some
problems in elementary geometry, in Hungarian), K\"oz. Mat. Lapok
61 (1980), 49--54.

\bibitem{E8} P.~Erd\H os and G.B.~Purdy,
\emph{Extremal problems in combinatorial geometry} in: Handbook of
Combinatorics, (R.L.~Graham, M.~Gr\"otschel, and L.~Lov\'asz eds.)
 Elsevier Science, (1995) 809--875.

\bibitem{Fal} G.~Faltings,
\emph{Endlichkeitss\"atze f\"ur abelsche Variet\"aten \"uber Zahlk\"orpern}
(Finiteness theorems for abelian varieties over number fields),
Invent. Math. 73 (3) (1983), 349--366.


\bibitem{Guy} R.~Guy, \emph{Unsolved Problems in Number Theory}
Problem Books in Mathematics Subseries: Unsolved Problems in
Intuitive Mathematics , Vol.1, Springer, 3rd ed., (2004) XVIII,
438 p.

\bibitem{HKM} H.~Harborth, A.~Kemnitz, and M.~M\"oller. \emph{An upper
bound for the minimum
diameter of integral point sets}, Discrete \& Comput. Geom.,
(1993) 9(4):427--432.

\bibitem{Hu} G.B.~Huff, \emph{Diophantine problems in geometry and
elliptic ternary forms},
Duke Math. J. vol. 15 (1948) 443--453.

\bibitem{Ke} A.~Kemnitz. \emph{Punktmengen mit ganzzahligen Abst\"anden},
Habilitationsschrift, TU Braunschweig, 1988.

\bibitem{KK}
T.~Kreisel and S.~Kurz, \emph{There are integral heptagons, no
three points on a line, no four on a circle}, Discrete \&
Computational Geometry, Online first:  DOI 10.1007/s00454-007-9038-6

\bibitem{Pe} W.~D.~Peeples Jr. \emph{Elliptic curves and rational
distance sets},
Proc. Am. Math. Soc., (1954) 5: 29--33.

\bibitem{Sil} J.~Silverman \emph{The Arithmetic of Elliptic
Curves}, Springer-Verlag, 1986


\bibitem{So} J.~Solymosi. \emph{Note on integral distances},
Discrete \& Comput. Geom., (2003) 30(2) 337--342.

\end{thebibliography}
\end{document}